\documentclass[twoside]{article}
\usepackage{amsmath,amsthm,amssymb,amscd}
\theoremstyle{plain}
 
\newtheorem{lem}{Lemma}[section]
\newtheorem{cor}{Corollary}[section]

\newtheorem{thm}{Theorem}[section]

\newcommand{\bb}[1]{\mbox{$\mathbb{#1}$}}
\theoremstyle{definition}
\newtheorem{ex}{Example}[section]
\theoremstyle{remark}
\newtheorem{rem}{Remark}[section]

\title{Hirzebruch classes and motivic Chern classes for singular (complex)
algebraic varieties}
\author{Jean-Paul Brasselet \and J\"{o}rg Sch\"{u}rmann \and Shoji Yokura}
\date{}

\begin{document}

\maketitle
\bibliographystyle{plain}

\begin{abstract}
In this paper we study some new theories of characteristic homology classes
of singular complex algebraic varieties. First we introduce a natural
transformation $T_{y}: K_{0}(var/X) \to H_{*}(X)\otimes \bb{Q}[y]$ commuting
with proper pushdown, which generalizes the corresponding Hirzebruch
characteristic. Here $K_{0}(var/X)$ is the relative Grothendieck group of
complex algebraic varieties over $X$ as introduced and studied by Looijenga
and Bittner in relation to motivic integration. $T_{y}$ is a homology class
version of the motivic measure corresponding to a suitable specialization
of the well known Hodge polynomial. This transformation unifies
the Chern class transformation of MacPherson and Schwartz (for $y=-1$)
and the Todd class transformation in the singular Riemann-Roch theorem of
Baum-Fulton-MacPherson (for $y=0$). In fact, $T_{y}$ is the
composition of a generalized version of this Todd class transformation due to
Yokura, and a new motivic Chern class transformation $mC_{*}: K_{0}(var/X) \to
G_{0}(X)\otimes \bb{Z}[y]$, which generalizes the total $\lambda$-class
$\lambda_{y}T^{*}X$ of the cotangent bundle to singular spaces. Here $G_{0}(X)$
is the Grothendieck group of coherent sheaves, and the construction of
$mC_{*}$ is based on some results from the theory of algebraic mixed Hodge
modules due to M.Saito. In the final part of the paper we use the algebraic
cobordism theory $\Omega_{*}$ of Levine and Morel to lift $mC_{*}$ further up
to a natural transformation $mC'_{*}: K_{0}(var/X) \to
\Omega_{*}(X)\otimes_{\bb{L}*} \bb{Z}[y,y^{-1}]$, where the last theory is a
suitable universal oriented Borel-Moore weak homology theory. Moreover, all our
results can be extended to varieties over a base field $k$ embeddable
into $\bb{C}$. \end{abstract}

\section*{Introduction}
In this paper we study some new theories of characteristic homology classes
of a singular algebraic variety $X$ defined over a base field $k$ embeddable
into $\bb{C}$. Let $K_{0}(var/X)$ be the relative Grothendieck group of
algebraic varieties over $X$ as introduced and studied by Looijenga \cite{Lo}
and Bittner \cite{Bi} in relation to motivic integration.\\

$ $\\
$ $\\

We introduce characteristic homomology class transformations
on $K_{0}(var/X)$ filling in the vertical maps in the following commutative
diagram:
\begin{displaymath} \begin{CD} 
K_{0}(var/X) @= K_{0}(var/X) @= \\ 
@VV mC'_{*} V  @VV mC_{*} V \\
\Omega_{*}(X)\otimes_{\bb{L}*} \bb{Z}[y,y^{-1}] 
@> \omega_{*} >> G_{0}(X)\otimes \bb{Z}[y,y^{-1}] @>>> 
\end{CD} \end{displaymath}
\begin{displaymath} \begin{CD}  
@= K_{0}(var/X) @= K_{0}(var/X)\\
@.  @VV  T_{y} V  @VV T_{1} V \\
@> td_{(1+y)} >>
H_{*}(X)\otimes \bb{Q}[y,y^{-1}] @>y=1 >>
H_{*}(X)\otimes \bb{Q} \:,
\end{CD} \end{displaymath}
Here $\omega_{*}$ comes from the universal property of the
oriented Borel-Moore weak homology theory 
$\Omega_{*}(X)\otimes_{\bb{L}*} \bb{Z}[y,y^{-1}]$ \cite{LM,Loe}, and
$td_{(1+y)}$ is the generalization given by Yokura \cite{Y4}
of the Todd class transformation $td_{*}$ used in the singular Riemann-Roch
theorem of Baum-Fulton-MacPherson \cite{BFM,FM} (for Borel-Moore homology)
or Fulton \cite{Ful} (for Chow groups).\\

All our transformation commute with proper (or projective) pushdown 
(in the case of $mC'_{*}$), and are uniquely determined by this property
together with the following normalization conditions for $X$ smooth and
pure $d$-dimensional:
\begin{itemize}
\item $mC'_{*}([id_{X}])= 
\tilde{\lambda}_{y}(T^{*}X)\cdot [X]$, with $[X]$ the fundamental class
and $\tilde{\lambda}_{y}$ a suitable total $\lambda$-class for this
context (as defined in the last section).
\item $mC_{*}([id_{X}])= \sum_{i=0}^{d} \; [\Lambda^{i} T^{*}X]\cdot y^{i}
= \lambda_{y}([T^{*}X])\cap [{\cal O}_{X}]$,
with $\lambda_{y}$ the total $\lambda$-class.
\item $T_{y}([id_{X}])=\widetilde{td_{(y)}}(TX) \cap [X]$,
with $\widetilde{td_{(y)}}$ the corresponding Hirzebruch characteristic
class.
\item $T_{1}([id_{X}])= L^{*}(TX)\cap [X]$, with $L^{*}$ the total
Thom-Hirzebruch characteristic $L$-class.
\end{itemize}

Moreover, the transformation $T_{y}$ fits also into a commutative
diagram
\begin{displaymath} \begin{CD}
F(X) @<e<<  K_{0}(var/X) @> mC_{0} >> G_{0}(X) \\
@V c_{*} VV  @V T_{y} VV   @VV td_{*} V \\
H_{*}(X)\otimes\bb{Q} @< y=-1 << H_{*}(X)\otimes\bb{Q}[y] @> y=0 >>
H_{*}(X)\otimes\bb{Q} \:.
\end{CD} \end{displaymath}
with $c_{*}$ the Chern-Schwartz-MacPherson transformation
\cite{Schwa,M1,BrS,Ken} on the group $F(X)$ of algebraically constructible
functions. Here $mC_{0}$ is the degree zero component of our motivic
Chern class transformation $mC_{*}$, and $e$ is simply given by
\[e([f: Y\to X]):=f_{!}1_{Y} \in F(X)\:.\]
Note that the homomorphisms $e$ and $mC_{0}$ are surjective.
So $T_{y}$ unifies the Chern-Schwartz-MacPherson
transformation $c_{*}$ and the Todd transformation $td_{*}$ of
Baum-Fulton-MacPherson (for Borel-Moore homology) or Fulton (for Chow groups).\\

Fix an embedding $k\subset \bb{C}$. The construction of
$mC_{*}$ is based on some results from the theory of algebraic mixed Hodge
modules due to M.Saito \cite{Sai1}-\cite{Sai6}, which imply the existence of the two natural
transformations 
\begin{displaymath} \begin{CD}
K_{0}(var/X) @> mH >> K_{0}(MHM(X/k))
@> gr^{F}_{-*}DR >>
G_{0}(X)\otimes \bb{Z}[y,y^{-1}] \,
\end{CD} \end{displaymath}
whose composition is our motivic transformation $mC_{*}$. 
Here $K_{0}(MHM(X/k))$ is the Grothendieck group of the abelian category
$MHM(X/k)$ of mixed Hodge modules on $X$ (depending on the
embedding $k\subset \bb{C}$), and $mH$ is defined
by  \[mH([f: Y\to X]):= [f_{!}\bb{Q}^{H}_{Y}]\in  K_{0}(MHM(X/k)) \:.\] 
Here $\bb{Q}^{H}_{Y}$ is in
some sense the 'constant Hodge module'  on $Y$. Similarly, $gr^{F}_{-*}DR$
comes form a suitable filtered de Rham complex of the filtered holonomic
D-module underlying a mixed Hodge module.\\

In the last section we lift our transformation $mC_{*}$ further up
to the natural tranformation $mC'$. This is based on the result of
Levine and Morel \cite{LM,Loe} that $\omega_{*}$ is an isomorphism for $X$ smooth
and pure-dimensional, together with a presentation of $K_{0}(var/X)$
in terms of a blow-up relation for smooth spaces, which is due to
Looijenga \cite{Lo} and Bittner \cite{Bi}.\\

It is a pleasure to thank P.Aluffi for some discussions
on this subject. The paper \cite{A} was a strong motivation for our work,
which started with the papers \cite{Y1}-\cite{Y6} of the third author.
Some of these papers are partly motivated by the final remark of
MacPherson's survey article
\begin{itemize} 
\item \cite[p.326]{M2}: '... It remains to be seen wether there is a unified
theory of characteristic classes of singular varieties like the classical one
outlined above. ...'.
\end{itemize}
We hope that our results give some key to MacPherson's
question and answer  the following question or problem:
\begin{itemize} 
\item \cite[p.267]{Y4}: '... Is there a theory of characteristic homology classes
unifying the above three characteristic homology classes of possible
singular varieties? ...'.
\item \cite[p.1]{A}: '... There is a strong motivic feel to the theory of
Chern-Schwartz-Mac\-Pherson classes, although this does not seem to have yet
been congealed into a precise statement in the literature. ...'.
\end{itemize}

\tableofcontents

\section{The generalized Hirzebruch theorem}

First we recall the classical generalized Hirzebruch Riemann-Roch theorem
\cite{Hi} (compare with \cite{Y3,Y4}). Let $X$ be a smooth complex projective
variety and $E$ a holomorphic vector bundle over $X$. Then the
$\chi_{y}$-characteristic of $E$ is defined by \begin{align*}
\chi_{y}(X,E):= 
&\sum_{p=0}^{<\infty} \chi(X,E\otimes \Lambda^{p}T^{*}X)\cdot y^{p}\\
=&\sum_{p=0}^{<\infty} \left( \sum_{i=0}^{<\infty}
(-1)^{i}dim_{C}H^{i}(X,E\otimes \Lambda^{p}T^{*}X) \right)\cdot y^{p} \:,
\end{align*}
with $T^*X$ the holomorphic cotangent bundle of $X$. Then one gets
\begin{equation} \label{eq:gHHR} 
\chi_{y}(X,E)= \int_{X} \widetilde{td_{(y)}}(TX)\cdot ch_{(1+y)}(E) \cap [X]
\quad \in \bb{Q}[y], \tag{gHRR}
\end{equation}
\[\text{with} \quad ch_{(1+y)}(E):= \sum_{j=1}^{rk\;E} e^{\beta_{j}(1+y)}\]
\[\text{and} \quad \widetilde{td_{(y)}}(TX):= \prod_{i=1}^{dim X}
Q_{y}(\alpha_{i})  \:.\]
Here $\beta_{j}$ are the Chern roots of $E$ and $\alpha_{i}$ are the Chern
roots of the tangent bundle $TX$. Finally $Q_{y}(\alpha)$ is the normalized
power series 
\[Q_{y}(\alpha):= \frac{\alpha(1+y)}{1-e^{-\alpha(1+y)}} -\alpha y
\quad \in \bb{Q}[y][[\alpha]] \:. \]
So this power series $Q_{y}(\alpha)$ specializes to
\begin{displaymath}
Q_{y}(\alpha) = 
\begin{cases}
\:1+\alpha &\text{for $y=-1$,}\\
\:\frac{\alpha}{1-e^{-\alpha}} &\text{for $y=0$,}\\
\:\frac{\alpha}{\tanh \alpha} &\text{for $y=1$.}
\end{cases} \end{displaymath}

Therefore the modified Todd class $\widetilde{td_{(y)}}(TX)$ unifies the
following important characteristic cohomology classes of $TX$:
\begin{displaymath}
\widetilde{td_{(y)}}(TX) = 
\begin{cases}
\:c^{*}(TX) &\text{the total Chern class for $y=-1$,}\\
\:td^{*}(TX) &\text{the total Todd for $y=0$,}\\
\:L^{*}(TX) &\text{the total Thom-Hirzebruch L-class for $y=1$,}
\end{cases} \end{displaymath}
and the gRRH-theorem for the trivial bundle $E$ specializes to
the calculation of the following important invariants
with $\chi_{y}(X):=\chi_{y}(X,{\cal O}_{X})$:
\begin{displaymath}
\chi_{y}(X) = 
\begin{cases}
\: e(X) = \int_{X} c^{*}(TX)\cap [X] &\text{the top. Euler
characteristic for $y=-1$,}\\ 
\: \chi(X) = \int_{X} td^{*}(TX)\cap [X] &\text{the arithmetic genus for
$y=0$,}\\ 
\: \sigma(X) = \int_{X} L^{*}(TX)\cap [X]  &\text{the signature for $y=1$.}
\end{cases} \end{displaymath}

These three invariants and classes have been generalized to a singular complex
algebraic variety $X$ in the following way (where the invariants are only
defined for $X$ compact):

\[e(X) = \int_{X} c_{*}(X), \quad \text{with} \quad
c_{*}: F(X) \to H_{*}(X):=
\begin{cases}
A_{*}(X)\\
H^{BM}_{2*}(X,\bb{Z})
\end{cases}\]
the Chern-Schwartz-MacPherson transformation \cite{Schwa,M1,BrS,Ken}
from the abelian group $F(X)$
of complex algebraically constructible functions to homology,
where one can use Chow groups $A_{*}(\cdot)$ or Borel-Moore homology groups
$H^{BM}_{2*}(\cdot,\bb{Z})$. Then $c_{*}(X):=c_{*}(1_{X})$.
\[\chi(X) = \int_{X} td_{*}(X), \quad \text{with} \quad
td_{*}: G_{0}(X) \to H_{*}(X)\otimes \bb{Q} \]
the Todd transformation in the singular Riemann-Roch theorem of
Baum-Fulton-MacPherson \cite{BFM,FM} (for Borel-Moore homology) or Fulton \cite{Ful}
(for Chow groups). Here $G_{0}(X)$ is the Grothendieck group
of coherent sheaves.
Then $td_{*}(X):=td_{*}([{\cal O}_{X}])$,
with $[{\cal O}_{X}]$ the class of the structure sheaf. 
Finally for compact $X$ one also has
\[\sigma(X) = \int_{X} L_{*}(X), \quad \text{with} \quad
L_{*}: \Omega(X) \to H^{BM}_{2*}(X,\bb{Q}) \]
the homology L-class transformation of Cappell-Shaneson \cite{CS1}
(as reformulated by Yokura \cite{Y1}). Here $\Omega(X)$ is the abelian group
of cobordism classes of selfdual constructible complexes.
Then $L_{*}(X):=L_{*}([{\cal IC}_{X}])$,
with $[{\cal IC}_{X}]$ the class of the intersection cohomology complex
of Goresky-MacPherson \cite{GM}.\\

All these transformations commute with the corresponding pushdown
for proper maps (where all spaces are assumed to be compact in the case
of the L-class transformation). They are uniquely characterized
by this pushdown property and the normalization condition, that for $X$ smooth
and pure-dimensional one gets back the corresponding classes of $TX$:
\[c_{*}(X)=c^{*}(TX)\cap [X],\quad td_{*}(X)=Td^{*}(TX)\cap [X] \quad
\text{and} \quad L_{*}(X)=L^{*}(TX)\cap [X] \:.\]
Here the uniqueness result follows from resolution of singularities.
So all these theories have the same formalism, but they are defined 
on completely different theories! Nevertheless, it is natural to ask
for another theory of characteristic homology classes of 
singular complex algebraic varieties, which unifies the above characteristic
homology classes (as in \cite{M2,Y3,Y4}). 
Of course in the smooth case, this is done by the 
generalized Todd class $\widetilde{td_{(y)}}(TX) \cap [X]$ of the tangent
bundle. We now explain a partial solution to this question.

\section{Hirzebruch characteristic classes for singular varieties}

In the following we consider algebraic varieties over a base field $k$
embeddable into  $\bb{C}$, and all spaces are assumed to be reduced. Then the
homology group  $H_{*}(X)$ is either the Chow group or the Borel-Moore homology
group in case $k=\bb{C}$. Let
$K_{0}(var/X)$ be the relative Grothendieck group of algebraic
varieties over $X$, i.e. the quotient of the free abelian group of isomorphism
classes of algebraic morphisms $Y\to X$ to $X$,  modulo the additivity
relation generated by  \begin{equation} \label{eq:add}
[Y\to X] = [Z\to Y \to X] + [Y\backslash Z \to Y \to X] 
\tag{add}
\end{equation}
for $Z\subset Y$ a closed algebraic subset of $Y$.
By resolution of singularities, $K_{0}(var/X)$ is generated by
classes $[Y\to X]$ with $Y$ smooth pure dimensional and proper
(or projective) over $X$! These relative groups were introduced by Looijenga
in his Bourbaki talk \cite{Lo} about motivic measures and motivic integration,
and then further studied by Bittner \cite{Bi}. From our point of view,
these are the 'motivic versions' of the group $F(X)$ of algebraically constructible
functions. In particular, they have the same formalism, i.e.
functorial pushdown $f_{!}$ and pullback $f^{*}$ for any algebraic map
$f:X'\to X$ (which is not necessarily proper), together with a ring
multiplication (with unit $[id_{X}]=c^{*}[id_{pt}]$ for $c:X\to pt$ a constant
map) satisfying the projection formula 
\[f_{!}(\alpha\cdot f^{*}\beta)=(f_{!}\alpha)\cdot \beta\]
and the base change formula $g^{*}f_{!}=f'_{!}g'^{*}$ for any cartesian
diagram
\begin{displaymath} \begin{CD}
Y' @>g'>> X' \\
@V f' VV  @VV f V \\
Y @> g >> X \:. 
\end{CD} \end{displaymath}
For later use, let us recall the simple definition of the pullback and
pushdown  for $f: X'\to X$, and of exterior products:
\[f_{!}([h:Z\to X'])=[f\circ h: Z\to X] \quad \text{and}\]
\[[Z\to X]\times [Z'\to X'] = [Z\times Z' \to X \times X']\:.\]
\[\text{Moreover,}\: f^{*}([g:Y\to X])= [g': Y'\to X']\]
is defined by taking fiber products as above.

By these exterior products, $K_{0}(var/pt)$ becomes a commutative ring and
$K_{0}(var/X)$ a $K_{0}(var/pt)$-module such that $f_{!}$ and $f^{*}$ are
$K_{0}(var/pt)$-linear.
 
\begin{thm} \label{thm:Ty}
There exists a unique group homomorphism $T_{y}$ commuting with pushdown
for proper maps:
\[T_{y}: K_{0}(var/X)\to H_{*}(X)\otimes\bb{Q}[y] \:,\]
satisfying the normalization condition
\[T_{y}([id_{X}])=\widetilde{td_{(y)}}(TX) \cap [X]\]
for $X$ smooth and pure-dimensional.
\end{thm}

Note also, that all our groups have a natural increasing filtration preserved
under  proper pushdown: $F_{k}K_{0}(var/X)$ is generated by $[X'\to X]$ with
$dim X' \leq k$. Similarly $F_{k}F(X)$ and $F_{k}G_{0}(X)$ is generated by 
constructible functions or classes of coherent sheaves with support of 
dimension $\leq k$. Finally we define
$F_{k}\left(H_{*}(X)\otimes\bb{Q}[y]\right):=
F_{k}H_{*}(X)\otimes F_{k}\bb{Q}[y]$, where each factor has its canonical
filtration coming from the natural grading. In particular, any evaluation
homomorphism $H_{*}(X)\otimes\bb{Q}[y]\to H_{*}(X)$ for an $y\in \bb{Q}$
is then also filtration preserving.

\begin{cor} \label{cor:Ty}
\begin{enumerate}
\item $T_{y}$ is filtration preserving.
\item $T_{y}$ commutes with exterior products.
\item One has the following Verdier Riemann-Roch formula for $f:X'\to X$
a smooth morphism (of constant relative dimension):
\[\widetilde{td_{(y)}}(T_{f}) \cap f^{*}T_{y}([Z\to X]) 
= T_{y}f^{*}([Z\to X]) \:.\]
Here $T_{f}$ is the bundle on $X'$ of tangent spaces to the fibers of $f$,
i.e. the kernel of the surjection $df: TX'\to f^{*}TX$.
In particular $T_{y}$ commutes with pullback under \'{e}tale morphisms
(i.e. smooth of relative dimension $0$).
\end{enumerate}
\end{cor}

\begin{proof} \begin{enumerate}
\item follows by induction on $dim X$ from resolution of
singularities, 'additivity' and the normalization condition.
\item follows  from the  normalization condition together with
\[(f\times f')_{*} \simeq f_{*}\times f'_{*}
\quad \text{on $K_{0}(var/X\times X')$ and $H_{*}(X\times X')
\otimes\bb{Q}[y]$}\] 
for $f,f'$ proper, and by the multiplicativity of
$\widetilde{td_{(y)}}$: 
\[\widetilde{td_{(y)}}(T(X\times X')) = 
\widetilde{td_{(y)}}(T(X)) \times \widetilde{td_{(y)}}(T(X')) \:.\]
Compare also with \cite{Kw,KwY} for the case of Chern classes. 
\item It is enough to prove the claim for $g: Z\to X$ proper
with $Z$ smooth (and pure dimensional). Then it
follows from the   projection formula 
\[g'_{*}(\alpha\cdot g'^{*}\beta)=(g'_{*}\alpha)\cdot \beta\] 
for $g'$ proper and the base
change formula $f^{*}g_{*}=g'_{*}f'^{*}$ for the cartesian diagram
\begin{displaymath} \begin{CD}
Z' @>f'>> Z \\
@V g' VV  @VV g V \\
X' @> f >> X 
\end{CD} \end{displaymath}
with $g,g'$ proper and $f,f'$ smooth (of constant relative dimension).
Here these formulae also hold for the homology $H_{*}(\cdot)\otimes\bb{Q}[y]$:
\begin{align*}
&\widetilde{td_{(y)}}(T_{f}) \cap f^{*}T_{y}([Z\to X]) =
\widetilde{td_{(y)}}(T_{f}) \cap f^{*}g_{*}T_{y}([id_{Z}]) \\ 
=&\widetilde{td_{(y)}}(T_{f}) \cap g'_{*}f'^{*}T_{y}([id_{Z}]) 
= g'_{*}\left( g'^{*}\widetilde{td_{(y)}}(T_{f})\cap f'^{*}(
\widetilde{td_{(y)}}(TZ)\cap [Z])\right)\\
=& g'_{*}\left( \widetilde{td_{(y)}}(T_{f'})\cup f'^{*}
\widetilde{td_{(y)}}(TZ)\cap [Z']\right)
= g'_{*} \left( \widetilde{td_{(y)}}(TZ')\cap [Z'] \right)\\
=& g'_{*}T_{y}([id_{Z'}]) = T_{y}([Z'\to X']) =  
T_{y}f^{*}([Z\to X]) \:.
\end{align*} 
Of course, we also used the multiplicativity and functoriality
of $\widetilde{td_{(y)}}$.
Compare also with \cite{Y5} for the case of Chern classes.
\end{enumerate}
\end{proof}

\begin{thm}
\label{thm:emC} There exist unique group homomorphisms
\[e: K_{0}(var/X)\to F(X) \quad \text{and} \quad 
mC_{0}: K_{0}(var/X)\to G_{0}(X) \]
 commuting with pushdown for proper maps and 
satisfying the normalization condition
\[e([id_{X}])=1_{X} \quad \text{and} \quad 
mC_{0}([id_{X}])=[{\cal O}_{X}] \] 
for $X$ smooth and pure-dimensional.
\end{thm}

\begin{cor} \label{cor:amC}
\begin{enumerate}
\item $e$ and $mC_{0}$ are filtration preserving.
\item $e$ and $mC_{0}$ commute with exterior products
and pullback for smooth morphisms (of constant relative dimension).
\item The following diagram of natural transformations commutes:
\end{enumerate}
\end{cor}
\begin{displaymath} \begin{CD}
F(X) @<e<<  K_{0}(var/X) @> mC_{0} >> G_{0}(X) \\
@V c_{*} VV  @V T_{y} VV   @VV td_{*} V \\
H_{*}(X)\otimes\bb{Q} @< y=-1 << H_{*}(X)\otimes\bb{Q}[y] @> y=0 >>
H_{*}(X)\otimes\bb{Q} \:.\quad \Box
\end{CD} \end{displaymath}
 
Note that the homomorphisms $e$ and $mC_{0}$ are surjective.
So by the diagram above, $T_{y}$ unifies the Chern-Schwartz-MacPherson
transformation $c_{*}$ \cite{Schwa,M1,BrS,Ken} and the Todd transformation $td_{*}$ of
Baum-Fulton-MacPherson  \cite{BFM,FM} (for Borel-Moore homology) or 
Fulton \cite{Ful} (for Chow groups).
Of course it is natural to use the induced Chern class transformation
with rational coefficients.

Here the transformation $e$ is simply given by
\[e([f: Y\to X]):= f_{!}(1_{Y})=f_{!}c^{*}1_{pt} \in F(X)\]
(with $c:X\to pt$ a constant map). In particular $e([id_{X}])=1_{X}$
also for any singular $X$ so that 
\[T_{-1}([id_{X}]) = c_{*}(X) \quad \text{for any singular $X$.}\]
So the corresponding Hirzebruch characteristic (class) of $X$ 
\begin{equation}
\chi_{y}(X):= T_{y}([X\to pt]) \quad \text{and} \quad
T_{y}(X):=T_{y}([id_{X}])  \end{equation}
is a natural lift of the Euler characteristic $e(X)$ (for $X$ complete)
and the Chern class $c_{*}(X)$ of $X$.\\

It seems that our Hirzebruch characteristic class 
$T_{y}(X)$ corresponds (for $k=\bb{C}$) to a similar class announced 
some years ago by Cappell and Shaneson \cite{CS2,Sh}. If this is the case,
then there is a mistake in their announcement, because they claim that
$T_{0}([id_{X}])$ is the singular Todd class $td_{*}(X)$ for any singular
$X$. But this is not (!) true for our class, i.e. $mC_{0}([id_{X}])\neq
[{\cal O}_{X}]\in G_{0}(X)$ for some singular spaces $X$.

\begin{ex} \label{ex:T0-Todd}
Let $X$ be a singular curve (i.e. $dim X=1$) such that $X$ is not maximal
(maximal is sometimes also called weakly normal), but the weak normalization
$X_{max}$ is smooth. Then the canonical projection $\pi : X':=X_{max}\to X$
is not an isomorphism, but nevertheless a topological homeomorphism.
By 'additivity' one gets
\[\pi_{*}([id_{X'}]) = [id_{X}] \in K_{0}(var/X) \quad\text{, so that}\]
\[T_{0}(X)=\pi_{*}T_{0}(X')=\pi_{*}td_{*}(X')=td_{*}\pi_{*}([{\cal O}_{X'}])
\:.\]
But by assumption $\pi_{*}([{\cal O}_{X'}])=[{\cal O}_{X}]+ 
n\cdot [{\cal O}_{pt}]$ with $n>0$ so that
\[T_{0}(X)= td_{*}(X)+ n\cdot [pt] \neq td_{*}(X) \in H_{*}(X)\otimes \bb{Q}
\:.\]
One gets similar examples in any dimension by taking the product with a
projective space.
\end{ex}

Taking a complete singular curve $X$ over $k=\bb{C}$ such that the
normalization $\pi: X':=X_{nor}\to X$ is not a topological homeomorphism, one
gets in the same way examples of singular $X$ with $T_{1}(X)\neq L_{*}(X)$.
Note that the normalization map $\pi$ is a small resolution of singularities
so that $\pi_{*}({\cal IC}_{X'})={\cal IC}_{X}$. So to distinguish between
these characteristic classes, we call (for later reasons) $T_{0}(X)$
the Hodge-Todd class and $T_{1}(X)$ the Hodge L-class of $X$.

\begin{ex} \label{ex:T0=Todd} 
Assume $k=\bb{C}$ and that the complex algebraic variety $X$ has
at most 'Du Bois singularities' in the sense of M.Saito \cite{Sai5}. 
For example $X$ has
only 'rational singularities', e.g. $X$ is a toric variety. Then
$mC_{0}([id_{X}]) = [{\cal O}_{X}]\in G_{0}(X)$ and therefore
$T_{0}(X)=td_{*}(X)$ (compare with the end of section \ref{sec:Hodge}).
\end{ex}

Using the 'additivity' in $K_{0}(var/X)$ and the natural transformation
$T_{y}$, one gets similar additivity properties of the Todd class
$td_{*}(X)=T_{0}(X)$ for $X$ smooth (or with a most 'Du Bois singularities'),
which seem to be new and do not follow directly from the original
definition.

\begin{ex} \label{ex:addTodd}
\begin{enumerate}
\item By the gHRR-theorem (which we recover in the algebraic context
in the next section), we have
$T_{y}([\bb{P}^{1}\to pt])=\chi_{y}(\bb{P}^{1})=1-y$ so that 
\[\chi_{y}(\bb{A}^{1})=-y \quad \text{and} \quad
\chi_{y}(\bb{P}^{n})= 1 -y + \dots + (-y)^{n} \]
by 'additivity' and 'multiplicativity' for exterior products.
\item $T_{y}$ becomes multiplicative in Zariski locally trivial bundles,
e.g. if $E\to X$ is an algebraic vector bundle of rank $r+1$, then
the corresponding projective bundle $\bb{P}(E)\to X$ is 
Zariski locally trivial so that
\[T_{y}([\bb{P}(E)\to X])= T_{y}([id_{X}]) \cdot (1 -y + \dots + (-y)^{r})
\in  H_{*}(X)\otimes \bb{Q}[y] \:.\]
\item Let $\pi: X'\to X$ be the blow-up of an algebraic variety $X$ along
an algebraic subvariety $Y$ such that the inclusion $Y\to X$ is a regular
embedding of pure codimension $r+1$
(e.g. $X$ and $Y$ are smooth). Then $\pi$ is an isomorphism over
$X\backslash Y$ and a projective bundle over $Y$ corresponding to the normal
bundle of $Y$ in $X$ of rank $r+1$. So by 'additivity' one gets 
\[T_{y}([X'\to X])= T_{y}([id_{X}])+ T_{y}([Y\to X]) \cdot 
(-y + \dots + (-y)^{r}) \in  H_{*}(X)\otimes \bb{Q}[y] \:.\]  
In particular $T_{0}([X'\to X])= T_{0}([id_{X}])$, which is a
homology class version of the birational invariance of the
arithmetic genus $\chi_{0}$. More generally, by pushing down to a
point we get for $X$ complete the blow-up formula
\begin{equation} \label{eq:blup}
\chi_{y}(X')= \chi_{y}(X) + \chi_{y}(Y) \cdot 
(-y + \dots + (-y)^{r}) \in  \bb{Q}[y] \:.
\end{equation}
\end{enumerate} \end{ex}

Note that in case 3. one also has $\pi_{*}td_{*}(X')=td_{*}(X)$,
by functoriality of $td_{*}$ and the relation 
$[R\pi_{*}\pi^{*}{\cal O}_{X}]=[{\cal O}_{X}]$ for such a blow-up.
Using the 'weak factorization theorem' \cite{AKMW,W}, we get the following result,
which seems to be new and was motivated by a corresponding study
of Aluffi about Chern classes \cite{A}.

\begin{cor} \label{cor:tdbirat}
Let $\pi: Y\to X$ be a resolution of singularities. Then the class 
\[\pi_{*}(Td^{*}(TY)\cap [Y])=\pi_{*}T_{0}([id_{Y}]) \in H_{*}(X)\otimes
\bb{Q}\]
is independent of $Y$.
\end{cor}

\begin{proof} Let $\pi: Y\to X$ and $\pi': Y'\to X$ be two
resolution of singularities, together with a resolution
of singularities of the fiber-product
$Z\to Y\times _{X} Y'$ so that we get induced birational
morphisms $p: Z\to Y$ and $p': Z\to Y'$.
By the 'weak factorization theorem' \cite{AKMW,W}, this map $p$ (or $p'$) can be
decomposed as a finite sequence of projections from smooth spaces lying over
$Y$ (or $Y'$), which are obtained by  blowing up or blowing down along  
smooth centers. By the birational invariance above we get
$\pi_{*}td_{*}(Z)= td_{*}(X)$ (or $\pi'_{*}td_{*}(Z)= td_{*}(X')$),
from which the claim follows.
\end{proof}

\begin{rem} \label{rem:local}
In the framework of motivic integration \cite{DL,Lo}, it is natural to localize
the $K_{0}(var/pt)$-module $K_{0}(var/X)$ with respect to the class
$\bb{L}:=[\bb{A}^{1}\to pt]$ of the affine line. Here the module structure
comes from pullback along a constant map $X\to pt$. Then $T_{y}$ induces a
a similar transformation
\[T_{y}: M(var/X):=K_{0}(var/X)[\bb{L}^{-1}]\to H_{*}(X)\otimes \bb{Q}[y,y^{-1}]\:,\]
since $\chi_{y}(\bb{A}^{1})=-y \in H_{*}(pt)\otimes \bb{Q}[y,y^{-1}]
= \bb{Q}[y,y^{-1}]$ is invertible. Note that the original transformation
$T_{y}$ is a ring homomorphism on a point space $pt$,
and a corresponding module homomorphism over any space $X$, by the
multiplicativity with respect to exterior products.
Similarly, $T_{y}$ extends to a transformation 
\[T_{y}: \widehat{M}(var/X) \to
 H_{*}(X)\otimes \bb{Q}[y][[y^{-1}]] \]
of the corresponding completions with respect to $'dim\to -\infty'$.
Here 
\[dim ([X'\to X]\bb{L}^{-n}) := dim(X')-n \]
and the completion $\widehat{M}(var/X)$ also comes up in the
context of motivic integration. 

In the absolute case $X=pt$
is was introduced by Kontsevich in his study of the 'arc-space'
${\cal L}(X)$ of $X$ as the value group of a 'motivic measure' 
$\tilde{\mu}$ on a suitable Boolean algebra of subsets of ${\cal L}(X)$.
This allows one (compare \cite[4.4]{DL})
 to introduce new invariants for $X$ pure-dimensional, but maybe
singular, as the value of 
\[\tilde{\mu}({\cal L}(X)) \in \widehat{M}(var/pt) \to R\]
under a suitable homomorphism to a ring $R$.
Instead of $\tilde{\mu}({\cal L}(X))$, one can also use related
'motivic integrals' over ${\cal L}(X)$.

By our work, one can now introduce a similar characteristic class as
\begin{displaymath} \begin{CD}
\tilde{\mu}/X\; ({\cal L}(X)) \in \widehat{M}(var/X) 
@> T_{y} >> H_{*}(X)\otimes \bb{Q}[y][[y^{-1}]]
\end{CD} \end{displaymath}
by using a 'relative motivic measure' $\tilde{\mu}/X$ with values in
$\widehat{M}(var/X)$ \cite[sec.4]{Lo}, and the same for 'motivic integrals'
(compare also with \cite{Y6}).  
\end{rem}
 
\section{Motivic Chern classes for singular varieties}

In this section we explain that our transformation $T_{y}$ is induced
by the generalized singular Riemann-Roch theorem of Yokura \cite{Y4} from another
transformation $mC_{*}$, which we call the 'motivic Chern class
transformation'.

\begin{thm} \label{thm:mC}
There exists a unique group homomorphism $mC_{*}$ commuting with pushdown
for proper maps:
\[mC_{*}: K_{0}(var/X)\to G_{0}(X)\otimes\bb{Z}[y,y^{-1}] \:,\]
satisfying the normalization condition
\[mC_{*}([id_{X}])= \sum_{i=0}^{dim X} \; [\Lambda^{i} T^{*}X]\cdot y^{i}
= \lambda_{y}([T^{*}X])\cap [{\cal O}_{X}]\]
for $X$ smooth and pure-dimensional.
Here $\lambda_{y}: K^{0}(X)\to K^{0}(X)\otimes\bb{Z}[[y]]$ is the total
$\lambda$-class transformation on the Grothendieck 
$K^{0}(X)$ of coherent locally free sheaves on $X$, with
$\cap [{\cal O}_{X}]: K^{0}(X)\to G_{0}(X)$ induced by 
$\otimes {\cal O}_{X}$, which is an isomorphism for $X$ smooth.
\end{thm}

\begin{cor} \label{cor:mC}
\begin{enumerate}
\item $mC_{*}$ is filtration preserving, if $G_{0}(X)\otimes\bb{Z}[y,y^{-1}]$
has the induced filtration coming from the grading with $y$ of degree one.
Moreover $mC_{*}$ maps to $G_{0}(X)\otimes\bb{Z}[y]$, and by projecting 
further on the component in degree $0$, we get the
transformation $mC_{0}$ of theorem \ref{thm:emC}.
\item $mC_{*}$ commutes with exterior products.
\item One has the following Verdier Riemann-Roch formula for $f:X'\to X$
a smooth morphism (of constant relative dimension):
\[\lambda_{y}(T^{*}_{f}) \cap f^{*}mC_{*}([Z\to X]) 
= mC_{*}f^{*}([Z\to X]) \:.\]
Here $f^{*}: G_{0}(X)\otimes\bb{Z}[y,y^{-1}] \to
G_{0}(X')\otimes\bb{Z}[y,y^{-1}]$
is induced from the corresponding pullback of Grothendieck groups
by linear extension over $\bb{Z}[y,y^{-1}]$.
In particular $mC_{*}$ commutes with pullback under \'{e}tale morphisms.
$\quad \Box$
\end{enumerate} \end{cor}

\begin{rem} \label{rem:local2}
The remarks of remark \ref{rem:local} apply in the same way also to
the motivic Chern class transformation $mC_{*}$.
\end{rem}

Let us now explain an important reformulation by Yokura \cite{Y4} of the singular
Riemann-Roch theorem of Baum-Fulton-MacPherson (for Borel-Moore homology)
and Fulton (for Chow groups). Consider the group homomorphism
\begin{equation} \label{eq:gBFM}
\begin{split}
&td_{(q)}: G_{0}(X)\otimes\bb{Z}[q,q^{-1}] \to H_{*}(X)\otimes\bb{Q}[q,q^{-1}]
\:,\\
&td_{(q)}([{\cal F}]):= \sum_{i=0}^{<\infty}\: td_{i}([{\cal F}])\cdot q^{-i}
\:, \end{split} \end{equation}
with $td_{i}$ the degree $i$ component of the transformation $td_{*}$,
which is linearly extended over $\bb{Z}[q,q^{-1}]$. Since $td_{*}$ is
degree preserving, this new transformation also commutes with proper pushdown
(which again is defined by linear extension over $\bb{Z}[q,q^{-1}]$).
Finally make the substitution $q=1+y$ to get the natural transformation
\begin{equation} \label{eq:gBFM2}
td_{(1+y)}: G_{0}(X)\otimes\bb{Z}[y,y^{-1}] \to
H_{*}(X)\otimes\bb{Q}[y,y^{-1}] \tag{gBFM}
\end{equation}
commuting with proper pushdown. By \cite{Y3} we have the important

\begin{lem} \label{lem:Yok}
Assume $X$ is smooth and pure $d$-dimensional. Then
\[td_{(1+y)}\left( \lambda_{y}(T^{*}X) \right) = \widetilde{td_{(y)}}(TX)
\cap [X] \in H_{*}(X)\otimes\bb{Q}[y]\:.\]
\end{lem}

Let us sketch the simple proof (since Yokura \cite[lem.(2.3.7)]{Y3} uses a different notation,
and works only in (Borel-Moore) homology for compact spaces over $\bb{C}$).
Since $X$ is smooth, the singular Riemann-Roch theorem reduces to the classical
Grothendieck Rie\-mann-Roch theorem, i.e. 
\begin{equation} \label{GRR}
td_{*}(E)= ch^{*}(E)\cup Td^{*}(TX)\cap [X]: G_{0}(X)=K^{0}(X) \to
H_{*}(X)\otimes \bb{Q}\:, \end{equation}
with $ch^{*}: K^{0}(X) \to H^{*}(X)\otimes \bb{Q}$ the Chern character.
So 
\[td_{i}(E) = \left( ch^{*}(E)\cup Td^{*}(TX) \right)^{d-i} \cap [X]\]
 and
\[\sum_{i=0}^{<\infty}\: td_{i}\left(  \lambda_{y}(T^{*}X)  \right)\cdot
(1+y)^{-i} = (1+y)^{-d}\cdot ch_{(1+y)} \left(  \lambda_{y}(T^{*}X)  \right)
Td_{(1+y)}(TX) \cap [X] \:.\]
If $\alpha_{i}$ are the Chern roots of $TX$, then the claim follows from
\begin{align*}
& (1+y)^{-d}\cdot ch_{(1+y)} \left(  \lambda_{y}(T^{*}X)  \right)
Td_{(1+y)}(TX) \\
&= (1+y)^{-d} \;\prod_{i=1}^{d}\;\left(1+ye^{-\alpha_{i}(1+y)}\right)
\cdot \frac{\alpha_{i}(1+y)}{1-e^{-\alpha_{i}(1+y)}} \\
&= \prod_{i=1}^{d}\;\left(\frac{\alpha_{i}(1+y)}{1-e^{-\alpha_{i}(1+y)}} 
-\alpha_{i}y \right) = \widetilde{td_{(y)}}(TX) \:. \quad \Box
\end{align*}

\begin{cor} \label{cor:TY=tdmC}
The natural transformation $T_{y}$ of theorem \ref{thm:Ty} is given as the 
composition
\[T_{y}:= td_{(1+y)} \circ mC_{*}: K_{0}(var/X)\to H_{*}(X)\otimes\bb{Q}[y]
\subset H_{*}(X)\otimes\bb{Q}[y,y^{-1}] \:.\]
So the motivic Chern class transformation $mC_{*}$ is a 'natural lift'
of the $T_{y}$ transformation. $\quad \Box$
\end{cor}

\begin{rem} \label{rem:not}
Maybe here is the right place to explain our notion 'motivic Chern class
transformation' for $mC_{*}$ (and $mC'_{*}$). Of course the notion
'motivic (dual) $\lambda$-class transformation' would also be possible
by the corresponding normalization condition for $X$ smooth.
But we understand our transformations $T_{y},mC_{*}$ and $mC'_{*}$
as natural 'motivic liftings' of the Chern-Schwartz-MacPherson transformation
$c_{*}$ by the following commutative diagram, with $T_{-1}([id_{X}])=c_{*}(X)$
for any singular $X$:
\end{rem}
\begin{displaymath} \begin{CD} 
K_{0}(var/X) @= K_{0}(var/X) @= \\ 
@VV mC'_{*} V  @VV mC_{*} V \\
\Omega_{*}(X)\otimes_{\bb{L}*} \bb{Z}[y,y^{-1}] 
@> \omega_{*} >> G_{0}(X)\otimes \bb{Z}[y,y^{-1}] @>>> 
\end{CD} \end{displaymath}
\begin{displaymath} \begin{CD}  
@= K_{0}(var/X) @> e >> F(X)\\
@.  @VV  T_{y} V  @VV c_{*} V \\
@> td_{(1+y)} >>
H_{*}(X)\otimes \bb{Q}[y,y^{-1}] @>y=-1 >>
H_{*}(X)\otimes \bb{Q} \:.
\end{CD} \end{displaymath}

Assume $X$ is smooth, pure $d$-dimensional and complete so that the constant
map $f: X\to pt$ is proper. Then (the proof of) lemma \ref{lem:Yok}, together
with the multiplicativity of the Chern character and the functoriality of
$td_{(1+y)}$ imply the classical (gHRR):
\begin{align*}
\chi_{y}(X,E) &= \chi\left(X,E\otimes \lambda_{y}(T^{*}X)\right) \\
& = td_{(1+y)}f_{!}\left(E\otimes \lambda_{y}(T^{*}X)\right) =
f_{!}\;td_{(1+y)}\left(E\otimes \lambda_{y}(T^{*}X)\right) \\
& = \int_{X}\: (1+y)^{-d}\cdot ch_{(1+y)} \left( E\otimes \lambda_{y}(T^{*}X) 
\right) Td_{(1+y)}(TX) \cap [X]\\
& = \int_{X}\: (1+y)^{-d}\cdot ch_{(1+y)} \left(\lambda_{y}(T^{*}X) 
\right) Td_{(1+y)}(TX) \cdot ch_{(1+y)}(E) \cap [X]\\
&= \int_{X} \widetilde{td_{(y)}}(TX)\cdot ch_{(1+y)}(E) \cap [X] \:.
\end{align*}

\section{Hodge theoretic definition of motivic Chern classes}
\label{sec:Hodge}

In this section we construct the motivic Chern class transformation 
$mC_{*}$ and therefore by the discussion before also the transformation
$T_{y}$ with the help of some fundamental results from the theory of
algebraic mixed Hodge modules due to M.Saito \cite{Sai1}-\cite{Sai6}.
Since this is a very complicated theory, we reduce our construction to a
few formal properties, together with a simple and instructive
explicite calculation for the normalization condition,
all of which are contained in the work of M.Saito.\\

Let us assume that our base field is $k=\bb{C}$.
To motivate the following constructions, let us first recall the
definition of the Hodge characteristic transformation
$Hc: K_{0}(var/pt)\to \bb{Z}[u,v] \:.$
By the now classical theory of Deligne \cite{De1,De2,Sr}, the cohomology groups
$V=H^{i}_{c}(X^{an},\bb{Q})$ of a complex algebraic variety
have a canonical  functorial mixed Hodge structure,
which includes in particular the following data on the finite dimensional
rational vector space $V$:
\begin{itemize}
\item A finite increasing (weight) filtration $W$ of $V$ with
$W_{i}=\{0\}$ for $i<<0$ and $W_{i}=V$ for $i>>0$.
\item A finite decreasing (Hodge) filtration $F$ of $V\otimes \bb{C}$ with
$F^{p}=V$ for $p<<0$ and $F^{p}=\{0\}$ for $p>>0$.
\end{itemize}
These filtrations have to satisfy some additional properties, which imply
that the transformation of taking suitable graded vector spaces
$gr^{W}_{i}, gr_{F}^{p}$ and $gr_{F}^{p}gr_{i}^{W}$ for $i,p\in \bb{Z}$
induce corresponding transformations on the Grothendieck group
$K^{0}(MHS)$ of the abelian category of (rational) mixed Hodge structures,
i.e. morphism of mixed Hodge structures are 'strictly stable' with respect
to the filtrations $F$ and $W$. Assume $Y$ is a closed algebraic subset of
$X$ with open complement $U:=X\backslash Y$. Then the maps in the long exact
cohomology sequence
\begin{equation} \label{eq:long}
\cdots \to H^{i}_{c}(U^{an},\bb{Q}) \to H^{i}_{c}(X^{an},\bb{Q}) \to
H^{i}_{c}(Y^{an},\bb{Q}) \to \cdots
\end{equation}
are morphisms of mixed Hodge structures so that the function
\begin{equation} \label{eq:addgr}
X\mapsto Hc(X):= \sum_{i,p,q=0}^{<\infty}\:(-1)^{i}(-1)^{p+q}\cdot 
dim_{C}\left(gr_{F}^{p}gr_{p+q}^{W}H^{i}_{c}(X^{an},\bb{C})\right) u^{p}v^{q}
\end{equation} 
satisfies the 'additivity' property (add). In this way we get the 
Hodge characteristic (compare \cite{Sr}):
\[Hc: K_{0}(var/pt)\to \bb{Z}[u,v],\: [X\to pt]\mapsto Hc(X).\]
Note that most references do not include the sign-factor
$(-1)^{p+q}$ in their definition of the Hodge characteristic! Our
sign convention fits better with the following normalization for
$X$ smooth and complete.
Specializing further, one gets also the (compare \cite{Lo})
\begin{itemize} \item
Hodge filtration characteristic corresponding to $(u,v)=(y,-1)$:
\[Hfc(X):=\sum_{i,p=0}^{<\infty}\:(-1)^{i} 
dim_{C}\left(gr_{F}^{p}H^{i}_{c}(X^{an},\bb{C})\right) (-y)^{p} \:.\] \item
Weight filtration characteristic corresponding to $(u,v)=(w,w)$: \[wc(X):=
\sum_{i,q=0}^{<\infty}\:(-1)^{i}\cdot 
dim_{C}\left(gr_{q}^{W}H^{i}_{c}(X^{an},\bb{Q})\right) (-w)^{q}\:.\] \item
Euler characteristic (with compact support)  corresponding to $(u,v)=(-1,-1)$:
\[e(X):= \sum_{i=0}^{<\infty}\:(-1)^{i}\cdot 
dim_{C}(H^{i}_{c}(X^{an},\bb{Q})) \:. \]
\end{itemize}
So these specializations fit into the following commutative diagram:
\begin{displaymath} \begin{CD}
\bb{Z}[u,v] @> u=y > v=-1 > \bb{Z}[y] \\
@V u=w V v=w V  @VV y=-1 V \\
\bb{Z}[w] @>> w=-1 > \bb{Z}\:.
\end{CD} \end{displaymath}
Finally, this classical Hodge theory \cite{De1,De2,Sr} implies for $X$ smooth complete
(of pure dimension $d$)
the 'purity result' $gr_{p+q}^{W}H^{i}(X^{an},\bb{C})=0$ for $p+q\neq i$,
together with
\begin{align*}
h^{p,q}(X):= &\sum_{i=0}^{<\infty}\: (-1)^{i}(-1)^{p+q}\cdot
dim_{C}\left(gr_{F}^{p}gr_{p+q}^{W}H^{i}_{c}(X^{an},\bb{C})\right)\\
= & dim_{C}\left(gr_{F}^{p}H^{p+q}(X^{an},\bb{C})\right)
= dim_{C}H^{q}(X^{an},\Lambda^{p}T^{*}X^{an})\\
= &  dim_{C}H^{q}(X,\Lambda^{p}T^{*}X)\:.
\end{align*}
Here the last two equalities follow from GAGA, the degeneration of the
'Hodge to de Rham spectral sequence'
\[E^{p,q}_{1}= H^{q}(X^{an},\Lambda^{p}T^{*}X^{an}) \to
H^{p+q}(X^{an},\Lambda^{\bullet}T^{*}X^{an})\]
at $E_{1}$ and the 'holomorphic Poincar\'{e} lemma' 
\[H^{p+q}(X^{an},\bb{C})\simeq H^{p+q}(X^{an},\Lambda^{\bullet}T^{*}X^{an})
\:.\] So the holomorphic de Rham complex 
$DR({\cal O}_{X^{an}}):=[\Lambda^{\bullet}T^{*}X^{an}]$ (with ${\cal
O}_{X^{an}}$ in degree zero) is a resolution of the constant sheaf $\bb{C}$ on
$X^{an}$, and the 'stupid decreasing filtration' 
\begin{equation} \label{eq:stup}
F^{p}DR({\cal O}_{X^{an}}):=[0 \to \cdots \to 0 \to \Lambda^{p}T^{*}X^{an} \to
\cdots  \Lambda^{d}T^{*}X^{an}] 
\end{equation}
induces the Hodge filtration $F$ on $H^{*}(X^{an},\bb{C})$.\\

In particular $T_{y}([X\to pt])= Hfc([X\to pt])$ for $X$ smooth and complete
by (gHRR). But these classes $[X\to pt]$ generate $K_{0}(var/pt)$ so
that we get the following Hodge theoretic description for any $X$:
\begin{equation} \label{eq:TyH}
\chi_{y}(X)= T_{y}([X\to pt]) = \sum_{i,p=0}^{<\infty}\:(-1)^{i} 
dim_{C}\left(gr_{F}^{p}H^{i}_{c}(X^{an},\bb{C})\right) (-y)^{p} \:.
\end{equation}

And exactly this description can be generalized to the context of
relative Grothendieck groups $K_{0}(var/X)$ using the machinery of
mixed Hodge modules of M.Saito. But before we explain this, let us
point out another remark. All our characteristics above are indeed
ring homomorphisms on $K_{0}(var/pt)$, because this is the
case for $\chi_{y}=Hfc$ (for example by remark \ref{rem:local}).
Such ring homomorphism are called 'characteristics' \cite{DL,Lo}
or sometimes also 'motivic measures' \cite{LL},
and there are much more examples known. In this sense
our transformation $T_{y}$ is certainly a 'motivic characteristic 
class' since it is a homology class version of the
motivic characteristic $Hfc$, just like the Chern-Schwartz-MacPherson
(class) transformation $c_{*}$ is the homology class version
of the Euler-Poincar\'{e} characteristic $e$. 

\begin{rem} \label{rem:ring}
Our 'motivic characteristic classes' are only group homomorphisms,
because the corresponding 'homology theories' are only groups for
a singular space $X$. But they commute with exterior products,
so that they are ring homomorphisms for $X=pt$.
But this is in general no longer true for $X$ smooth,
even if one has then a corresponding ring structure on the
(co)homology. This is closely related to a corresponding
Verdier Riemann-Roch formula for the diagonal embedding 
$d: X\to X\times X$ (compare \cite{Sch2,Y5}).
\end{rem}

One can ask if also the other characteristics
$Hc$ or $wc$ can be 'lifted up' to such a homology class
transformation. But here the answer will be no (compare \cite{Jo})!

\begin{ex} \label{ex:counter}
Assume that there is a functorial transformation 
\[T_{u,v}: K_{0}(var/X)\to H_{*}(X)\otimes\bb{Q}[u,v]\]
commuting with proper
pushdown, and also with pullback $f^{*}$ for a finite smooth morphism
$f: X' \to X$ between smooth varieties, such that for $X=pt$ we get back the
Hodge characteristic $Hc$. Let $d$ be the degree of such a covering map
$f$. Then it follows (with $c: X\to pt$ a constant map):
\[Hc(X')= T_{u,v}([X'\to pt]) = T_{u,v}(c_{*}f_{*}f^{*}[id_{X}]) =
c_{*}f_{*}f^{*}T_{u,v}([id_{X}])\]
\[ = c_{*}(d\cdot T_{u,v}([id_{X}])) = d \cdot c_{*}T_{u,v}([id_{X}])
= d\cdot Hc(X) = Hc(d\cdot [id_{pt}]) \cdot Hc(X)  \:.\]
So (as usual), the transformation $Hc$ has then to be multiplicative
in such finite coverings. But this is not the case.
Let $X'\to X$ be such a finite covering of degree $d>1$ over an elliptic
curve $X$. Then $X'$ is also an elliptic curve so that 
\[  Hc(X) = Hc(X') = (1+u)(1+v) \neq 0\:.\]
Note that the same argument applies also to the weight characteristic,
with $wc(X)=wc(X')=(1+w)^{2}\neq 0$. Of course, everything is ok
in the context of $e$ and $Hfc$, since both are zero for an elliptic 
curve!
\end{ex}

Let us now formulate those results about algebraic mixed Hodge modules,
which we need for our application to the motivic Chern class transformation
$mC_{*}$. All these results are contained in the deep and long work of M.Saito.
Since most readers will not be familiar with this theory,
we present them in an axiomatic way pointing out the similarities
to constructible functions $F(X)$ and motivic Grothendieck groups 
$K_{0}(var/X)$.\\

Let $k$ be a subfield of $\bb{C}$. Then we work in the category
of reduced seperated schemes of finite type over $spec(k)$,
which we also call 'spaces' or 'varieties', with $pt=spec(k)$.

\begin{enumerate}
\item[MHM1:] To such a space $X$ one can associate an abelian category
of (algebraic) mixed Hodge modules $MHM(X/k)$, together with a functorial
pullback $f^{*}$ and pushdown $f_{!}$ on the level of derived categories
$D^{b}MHM(X/k)$ for any (not necessarily proper) map \cite[sec.4]{Sai2}
(and compare also with \cite{Sai5,Sai6}).
These transformations are functors of triangulated categories.
\item[MHM2] Let $i:Y\to X$ be the inclusion of a closed subspace, 
with open complement $j: U:=X\backslash Y \to X$.
Then one has for $M\in D^{b}MHM(X/k)$ a distinguished triangle (\cite[(eq.(4.4.1),p.321]{Sai2})
\[j_{!}j^{*}M \to M \to i_{!}i^{*}M \stackrel{[1]}{\to} \:.\]
\item[MHM3:] For all $p\in \bb{Z}$ one has a functor of triangulated
categories
\[gr^{F}_{p}DR: D^{b}MHM(X/k) \to D^{b}_{coh}(X)\]
commuting with proper pushdown (compare with \cite[sec.2.3]{Sai1},
\cite[p.273]{Sai2}, \cite[eq.(1.3.4),p.9, prop.2.8]{Sai5} and also with 
\cite{Sai4,Sai5}). Here $D^{b}_{coh}(X)$ is
the bounded derived category of sheaves of ${\cal O}_{X}$-modules 
with coherent cohomology sheaves. Moreover, $gr^{F}_{p}DR(M)=0$
for almost all $p$ and $M\in D^{b}MHM(X/k)$ fixed (\cite[prop.2.2.10,
eq.(2.2.10.5)]{Sai1} and \cite[lem.1.14]{Sai3}).
\item[MHM4:] There is a distinguished element $\bb{Q}_{pt}^{H}
\in MHM(pt/k)$ such that 
\[gr^{F}_{-p}DR(\bb{Q}_{X}^{H}) \simeq \Lambda^{p}T^{*}X[-p]
\in D^{b}_{coh}(X)\]
for $X$ smooth and pure dimensional (\cite{Sai2}).
Here $\bb{Q}_{X}^{H}:=c^{*}\bb{Q}_{pt}^{H}$ for $c: X\to pt$
a constant map, with $\bb{Q}_{pt}^{H}$ viewed as a complex
concentrated in degree zero.
\end{enumerate}

Since the above transformations are functors of triangulated
categories, they induce functors on the level of Grothendieck groups
(of triangulated categories)
which we denote by the same name. By general nonsense one gets
for these Grothendieck groups isomorphisms
\[K_{0}(D^{b}MHM(X/k)) \simeq K_{0}(MHM(X/k)) \quad \text{and}
\quad K_{0}(D^{b}_{coh}(X)) \simeq G_{0}(X)\]
by associating to a complex its alternating sum of cohomology objects.

As explained later on, $K_{0}(MHM(X/k))$ plays the role of
'Hodge constructible functions' with the class of $\bb{Q}_{X}^{H}$
as 'constant Hodge function on $X$'.
By (MHM3) we get a group homomorphism commuting with proper pushdown:

\begin{equation} \label{eq:grH}
\begin{split}
 gr^{F}_{-*}DR: K_{0}(MHM(X/k)) \to G_{0}(X)\otimes \bb{Z}[y,y^{-1}] \:;\\
[M] \mapsto \sum_{p} \: [gr^{F}_{-p}DR(M)]\cdot (-y)^{p} \:.
\end{split} \end{equation}

And as for the map $e$ from motivic Grothendieck groups to constructible
functions, we get a group homomorphism commuting with pushdown
(compare also with \cite[sec.4]{Lo}):
\begin{equation} \label{eq:mH}
mH: K_{0}(var/X) \to K_{0}(MHM(X/k))\:, \:
[f: X'\to X] \mapsto [f_{!} \bb{Q}_{X'}^{H}]  \:.
\end{equation}

Indeed, the 'additivity relation' (add) follows from (MHM2)
together with the functoriality of pushdown and pullback (MHM1):
For $i: Y\to X'$ the inclusion of a closed subspace, with open
complement $j:U\to X'$, the distinguished triangle
\[j_{!}j^{*}c^{*}\bb{Q}_{pt}^{H} \to c^{*}\bb{Q}_{pt}^{H} 
\to i_{!}i^{*}c^{*}\bb{Q}_{pt}^{H} \stackrel{[1]}{\to} \]
induces under $f_{!}$ the distinguished triangle (with $f: X'\to X$ as before)
\[f_{!}j_{!}j^{*}c^{*}\bb{Q}_{pt}^{H}  \to f_{!}c^{*}\bb{Q}_{pt}^{H} 
\to f_{!}i_{!}i^{*}c^{*}\bb{Q}_{pt}^{H} \stackrel{[1]}{\to} \:.\]
It translates in the corresponding Grothendieck group into the relation
\[[f_{!}c^{*}\bb{Q}_{pt}^{H}] = [f_{!}j_{!}j^{*}c^{*}\bb{Q}_{pt}^{H}]
+ [f_{!}i_{!}i^{*}c^{*}\bb{Q}_{pt}^{H}]\:.\]
This finally is nothing else than the asked additivity property 
\[[f_{!} \bb{Q}_{X'}^{H}] = [(f\circ j)_{!}\bb{Q}_{U}^{H}] +
 [(f\circ i)_{!}\bb{Q}_{Y}^{H}] \in K_{0}(D^{b}MHM(X/k)) \:.\]
Moreover, $mH$ commutes with pushdown 
for a map $f: X'\to X$ again by functoriality:
\[mH(f_{!}[g:Y\to X']) = mH([f\circ g: Y\to X]) =
[(f\circ g)_{!}\bb{Q}_{Y}^{H}] \]
\[= [f_{!}g_{!}\bb{Q}_{Y}^{H}] = f_{!}[g_{!}\bb{Q}_{Y}^{H}] 
= f_{!}mH([g:Y\to X])\:.\]

By (MHM4) we get for $X$ smooth and pure dimensional:
\[ gr^{F}_{-*}DR\circ mH([id_{X}]) =  
\sum_{i=0}^{dim X} \; [\Lambda^{i} T^{*}X]\cdot y^{i}
\in G_{0}(X)\otimes \bb{Z}[y,y^{-1}] \:.\]

\begin{cor} \label{cor:mc=DRmH}
The motivic Chern class transformation $mC_{*}$ of theorem \ref{thm:mC}
is given as the composition $mC_{*}= gr^{F}_{-*}DR\circ mH$.
$\quad \Box$
\end{cor}

\begin{rem} \label{rem:emb}
The definition of $MHM(X/k)$ and therefore also the transformations
$mH$ and  $gr^{F}_{-*}DR$ depend a priory on the embedding $k\subset \bb{C}$.
By the uniqueness statement of Theorem \ref{thm:Ty} and \ref{thm:mC}
this is not the case for the transformations $T_{y}$ and $mC_{*}$,
i.e. they are independent of the choice of the embedding $k\subset \bb{C}$
used in their definition in Corollary \ref{cor:TY=tdmC} and \ref{cor:mc=DRmH}.
\end{rem}

Let us now explain a little bit of the definition of the abelian category
$MHM(X/k)$ of algebraic mixed Hodge modules on $X/k$.
Its objects are special tuples $((M,F),K,W)$, which for $X$ smooth are given by
\begin{itemize}
\item $(M,F)$ an algebraic holonomic filtered D-module $M$ on $X$ with an
exhaustive,  bounded from below and increasing (Hodge) filtration $F$ by
algebraic ${\cal O}_{X}$-modules such that $gr_{*}^{F}M$ is a coherent
$gr_{*}^{F}{\cal D}_{X}$-module. In particular, the filtration $F$ is
 finite, which will imply the last claim of (MHM3).  Here the filtration $F$ on
the sheaf of algebraic differential operators ${\cal D}_{X}$ on $X$ is the
order filtration, and one can work either with left or right D-modules.
For singular $X$ one works with suitable local embeddings into manifolds
and corresponding filtered D-modules with support on $X$
(compare \cite{Sai1,Sai2,Sai4}).
\item $K\in D^{b}_{const}(X(\bb{C})^{an},\bb{Q})$ is an algebraically constructible
complex of sheaves of $\bb{Q}$-vector spaces (with finite dimensional stalks,
compare for example with \cite{Sch1}))
on the associated analytic space $X(\bb{C})^{an}$ corresponding to the induced
algebraic variety $X(\bb{C}):=X\otimes_{k}\bb{C}$ over $\bb{C}$,
which is perverse with respect to the middle perversity $t$-structure.
$F$ is called the underlying rational sheaf complex.
\item In addition one fixes a quasi-isomorphism $\alpha$ between 
$K(\bb{C}):=K \otimes_{Q}\bb{C}$ and the holomorphic de Rham complex
$DR(M(\bb{C})^{an})$ associated to the induced ${\cal D}_{X(C)}$-module
$M(\bb{C}):=M \otimes_{k}\bb{C}$.
\item $W$ is finally a finite increasing (weight) filtration of
$(M,F)$ and $K$, compatible in the obvious sense with the quasi-isomorphism
$\alpha$ above.
\end{itemize}
These data have to satisfy a long list of properties which we do not recall
here (since it is not important for us). In particular, one gets the
equivalence \cite[eq.(4.2.12),p.319]{Sai2})
\begin{equation} \label{eq:MHM=MHS}
\begin{split}
MHM(pt/\bb{C}) &\simeq \{\text{(graded) polarizable mixed $\bb{Q}$-Hodge
structures} \} \\
\:F^{-p} &\leftrightarrow F_{p} \end{split}
\end{equation}
between the category of algebraic mixed Hodge modules on $pt=spec(\bb{C})$,
and the category of (graded) polarizable mixed $\bb{Q}$-Hodge
structures. Of course, one has to switch the increasing $D$-module
filtration $F^{p}$ into an decreasing Hodge filtration by
$F^{-p}\leftrightarrow F_{p}$  so that $gr^{F}_{-p}\simeq gr_{F}^{p}$.
For elements in $MHM(pt/k)$, the corresponding Hodge filtration is already
defined over $k$ (compare \cite[sec.1.3]{Sai6}). \\
The distinguished element $\bb{Q}^{H}_{pt}\in MHM(pt/k)$ of (MHM4) is given by
\begin{equation} \label{eq:disting}
\bb{Q}^{H}_{pt}:=((k,F),\bb{Q},W) \quad \text{with $gr^{F}_{i}=0=gr^{W}_{i}$
for $i\neq 0$}
\end{equation}
and $\alpha: k \otimes \bb{C}\simeq \bb{Q}\otimes \bb{C}$ the obvious
isomorphism (compare \cite[sec.1.3]{Sai6}).
The functorial pullback and pushdown of (MHM1) corresponds under the forget
functor (\cite[thm.0.1,p.222]{Sai2})
\begin{equation} \label{eq:rat}
rat: D^{b}MHM(X/k)\to D^{b}_{constr}(X(\bb{C})^{an},\bb{Q}) \:,
((M,F),K,W)\mapsto K 
\end{equation}
to the classical corresponding (derived) functors $f^{*}$ and $f_{!}$
on the level of algebraically constructible sheaf complexes,
with $rat(\bb{Q}^{H}_{X})\simeq \bb{Q}_{X(C)^{an}}$.
So by (\ref{eq:MHM=MHS}) one should think of an algebraic mixed Hodge module
as a kind of '(perverse) constructible Hodge sheaf'!
But one has to be very careful with this analogy. $\bb{Q}^{H}_{X}$
is in general a highly complicated complex in $D^{b}MHM(X/k)$,
which is impossible to calculate explicitely. But if $X$ is smooth and
pure $d$-dimensional, then $\bb{Q}_{X(C)^{an}}[-d]$ is a perverse sheaf
and $\bb{Q}^{H}_{X}[-d]\in MHM(X/k)$ a single mixed Hodge module
(in degree $0$), which is explicitly given by (\cite[eq.(4.4.2),p.322]{Sai2}):
\begin{equation} \label{MHMsmooth}
\bb{Q}^{H}_{X}[-d] \simeq (({\cal O}_{X},F),\bb{Q}_{X(C)^{an}}[-d],W) \:,
\end{equation}
with $F$ and $W$ the trivial filtration $gr_{i}^{F}=0=gr_{i}^{W}$
for $i\neq 0$. Here we use for the underlying D-module the description as a
left D-module, which maybe is more natural at this point.\\

The distingished triangle (MHM2) is a 'lift' of the corresponding
distinguished triangle for constructible sheaves.
Similarly, by taking a constant map $f: X\to pt$ we get by (MHM1) and
(\ref{eq:MHM=MHS}) a functorial (rational) mixed Hodge structure on
\[rat \left(R^{i}f_{!}\bb{Q}^{H}_{X} \right) \simeq
H^{i}_{c}(X(\bb{C})^{an},\bb{Q}) \:,\]
whose Hodge numbers are easily seen to be the same as those coming
from the mixed Hodge structure of Deligne \cite{De1,De2,Sr} (both have the same additivity
property so that one only has to compare the case $X$ smooth, which
follows from the constructions). In fact, even the Hodge structures 
are the same by a deep theorem of M.Saito \cite[cor.4.3]{Sai5}.\\

Let us finally explain (MHM3) and (MHM4) for the case $X$ smooth and
pure $d$-dimensional. The de Rham functor $DR$ factorizes as
(compare \cite{Sai1,Sai3,Sai5})
\begin{equation} \label{eq:facDR}
DR: D^{b}MHM(X/k) \to D^{b}F_{coh}(X,Diff) \to  
D^{b}_{const}(X(\bb{C})^{an},\bb{C}) \:,
\end{equation}
with $D^{b}F_{coh}(X,Diff)$ the 'bounded derived category of
filtered differential complexes' on $X$ with coherent cohomology sheaves.
Here the objects are bounded complexes $(L^{\bullet},F)$ of ${\cal
O}_{X}$-sheaves with an increasing (bounded from below) filtration $F$ by such
sheaves, whose morphisms are 'differential operators' in a suitable sense. In
particular \[gr^{F}_{p}(L^{\bullet}) \in D^{b}(X,{\cal O}_{X})\]
becomes an ${\cal O}_{X}$-linear complex with coherent cohomology.
Moreover, the morphisms of mixed Hodge modules are 'strict' with respect to the
Hodge filtration $F$ (and the weight filtration $W$) so that 
$gr^{F}_{p}DR$ induces the corresponding transformation of (MHM3).
Finally, $DR(\bb{Q}^{H}_{X})$ is given by the usual de Rham complex
$\Lambda^{\bullet}T^{*}X$ with the induced increasing filtration  
\begin{equation} \label{eq:DR2}
F_{p}DR(\bb{Q}^{H}_{X}) := [ F_{p}{\cal O}_{X} \to 
F_{p+1}{\cal O}_{X}\otimes \Lambda^{1} T^{*}X \to \: \cdots \to  
F_{p+d}{\cal O}_{X}\otimes \Lambda^{d} T^{*}X ]\:,
\end{equation}
with $F_{p}{\cal O}_{X}$ in degree zero and $F$ the trivial filtration with
$gr_{i}^{F}=0$ for $i\neq 0$. If we switch to the corresponding decreasing
filtration (with $gr^{F}_{-p}\simeq gr^{p}_{F}$):
\[F^{p}DR(\bb{Q}^{H}_{X}) := [F^{p}{\cal O}_{X} \to 
F^{p-1}{\cal O}_{X}\otimes \Lambda^{1} T^{*}X \to \: \cdots \to  
F^{p-d}{\cal O}_{X}\otimes \Lambda^{d} T^{*}X] \:,\]
then one gets to this complex a natural filtered map form the de Rham complex
$DR({\cal O}_{X})$ with the stupid filtration $\sigma^{p}$ as in
(\ref{eq:stup}): 
\[ [0 \to \cdots \to 0 \to \Lambda^{p} T^{*}X \to
\cdots  \Lambda^{d} T^{*}X] \:.\]
And one trivially checks that this induces on the associated graded
complexes the isomorphism 
\[gr^{p}_{\sigma} DR({\cal O}_{X}) \simeq \Lambda^{p}T^{*}X[-p] 
 \simeq gr^{p}_{F}DR (\bb{Q}^{H}_{X}) 
\simeq gr_{-p}^{F}DR (\bb{Q}^{H}_{X}) \:.\]
In this way one finally also gets (MHM4). 

\begin{rem} \label{rem:totaro}
The use of the transformation $gr^{F}_{p}DR$ of (MHM3) in the context
of characteristic classes of singular spaces is not new.
It was already used by Totaro \cite{To} in his study of the relation between
Chern numbers for singular complex varieties and elliptic homology!
But he was interested in characteristic numbers and classes
invariant under 'small resolution', and not in functoriality as in
our paper. So  he works with the counterpart ${\cal IC}_{X}^{H}
\in MHM(X/\bb{C})$ of the intersection cohomology complex
instead of the constant 'Hodge sheaf' 
$\bb{Q}_{X}^{H} \in D^{b}MHM(X/\bb{C})$ as used in this paper.
He then also applied the singular Riemann-Roch transformation
$td_{*}$ of Baum-Fulton-MacPherson to associate to a singular
complex algebraic variety $X$ of dimension $n$ some natural homology classes
$\chi^{n-k}_{p}(X)\in H^{BM}_{2k}(X,\bb{Q})$ for $p\in \bb{Z}$.
In our notation, the corresponding total homology class 
$\chi^{n-*}_{p}(X)\in H^{BM}_{2*}(X,\bb{Q})$ 
is given by evaluating
\[td_{(1+y)}\circ gr^{F}_{-p}DR({\cal IC}_{X}^{H}) \in 
H^{BM}_{2*}(X,\bb{Q})[y,y^{-1}]\]
at $y=0$. Here it is important to work with the
transformation 
\[td_{(1+y)}\circ gr^{F}_{-p}DR: K_{0}(MHM(X/\bb{C}))\to 
H^{BM}_{2*}(X,\bb{Q})[y,y^{-1}] \:.\]
This allows one to use more general
coefficients like ${\cal IC}_{X}^{H} \in MHM(X/\bb{C})$,
which are in general not in the image of $mH$.
\end{rem}

As explained before, it is in general impossible to calculate
$gr_{-p}^{F}DR (\bb{Q}^{H}_{X})$ explicitely for a singular
space $X$. But by the work of M.Saito \cite{Sai5}, one has for a
complex algebraic variety $X$ of dimension $n$
at least some general informations: 
\begin{enumerate}
\item $gr_{-p}^{F}DR (\bb{Q}^{H}_{X})\simeq 0 \in D^{b}_{coh}(X)$
for $p<0$ and $p>n$. Of course, the corresponding vanishing in
$G_{0}(X)$ follows also from Corollary \ref{cor:mc=DRmH}.
\item Let $X'$ be a resolution of singularities of the union
of the $n$-dimensional irreducible components of $X$, with
$\pi: X'\to X$ the induced proper map. Then
\[ gr_{-n}^{F}DR (\bb{Q}^{H}_{X})\simeq \pi_{*} 
\Lambda^{n}T^{*}X'[-n] \:.\]
Note that $R^{i}\pi_{*} \Lambda^{n}T^{*}X'=0$ for $i>0$ by the
Grauert-Riemenschneider vanishing theorem. 
\item $h^{i}(gr_{0}^{F}DR (\bb{Q}^{H}_{X}))\simeq 0$ for $i<0$,
and 
\[h^{0}(gr_{0}^{F}DR (\bb{Q}^{H}_{X}))\simeq {\cal O}_{X}^{wn} \:,\]
with ${\cal O}_{X}^{wn}$ the coherent structure sheaf 
of the weak normalization $X_{max}$ of $X$ (whose underlying space
is identified with $X$). One gets in particular natural morphisms
\[{\cal O}_{X} \to {\cal O}_{X}^{wn} \to
gr_{0}^{F}DR (\bb{Q}^{H}_{X})\]
in $D^{b}_{coh}(X)$. $X$ has by definition 'at most Du Bois singularities',
if these morphisms are quasi-isomorphisms.
This is by for example the case, if $X$ has at most 'rational singularities'
(\cite[thm.5.4]{Sai5}), i.e.
$X$ is normal and $R^{i}\pi_{*} {\cal O}_{X'}=0$ for $i>0$ for some
(and then for any) resolution of singularities $\pi: X'\to X$.
\end{enumerate}

\section{Motivic Chern classes in an algebraic cobordism theory}

We continue to work on the category of reduced seperated schemes
of finite type over a base field $k$ embeddable into $\bb{C}$.
Since we want to use some results from Levine and Morel \cite{LM,Loe} about
their algebraic cobordism theory $\Omega_{*}$, we use the convention
of that paper that all smooth morphisms of this section are assumed
to be quasi-projective! In particular, a smooth variety is assumed
to be quasi-projective.\\

Let $M^{+}_{*}(X)$ be the free abelian group of isomorphism classes
of projective morphism $f: X'\to X$ with $X'$ smooth and pure-dimensional,
graded by the dimension of $X'$. By Chow's lemma and resolution of
singularities one gets a group epimorphism
\[v: M^{+}_{*}(X)\to K_{0}(var/X)\] 
by associating to such an $f$ its class
$[f:X'\to X]$ in the Grothendieck group of algebraic varieties over $X$.
By a deep theorem of Bittner \cite{Bi}  and Looijenga (relying on the
weak factorization theorem \cite{AKMW,W}), one gets $K_{0}(var/X)$ form $M^{+}_{*}(X)$ by 
imposing the 'blow-up relation' 
\begin{equation} \label{eq:bl}
[\emptyset \to X]=0 \quad \text{and} \quad
[Bl_{Y}X'\to X] - [E\to X]= [X'\to X] - [Y\to X] \:,\tag{bl}
\end{equation}
with $X'\to X$ as before and $Y\subset X'$ a closed pure-dimensional smooth
subspace of $X'$. Here $Bl_{Y}(X')$ is the blow-up of $X'$ along $Y$ 
with exceptional divisor $E$. These are then also smooth pure-dimensional
and the induced maps to $X$ are projective, so that it makes sense
to impose the relation (bl) on $M^{+}_{*}(X)$.
Of course it is clear that the relation (bl) holds in 
$K_{0}(var/X)$, since $Bl_{Y}X'\backslash E$ is isomorphic to $X'\backslash Y$.
Note that this relation (bl)
is not comparable with the graduation of $M^{+}_{*}(X)$,
but with the corresponding filtration so that $v$ is a filtered 
map.\\

One also gets the algebraic cobordism group $\Omega_{*}(X)$ of $X$
out of $M^{+}_{*}(X)$ by imposing other kinds of relations \cite{LM,Loe}
 in such a way
that 
\[co: M^{+}_{*}(X) \to \Omega_{*}(X)\]
is a graded group epimorphism (compare \cite[lem.4.17,p.45]{LM})
and the theory $\Omega_{*}$ becomes an universal
oriented Borel-Moore weak homology theory in the sense of Levine-Morel.
Here an oriented Borel-Moore weak homology theory $A_{*}$ is given
by a transformation $X\mapsto A_{*}(X)$ with $A_{*}(X)$ a commutative graded
group such that one has:
\begin{itemize}
\item $A_{*}$ is additive in the sense that it maps finite disjoint
unions isomorphic into finite products.
\item $A_{*}$ is covariant functorial with a pushdown $f_{*}$ for $f$
projective, and contravariant functorial with a pullback $f^{*}$
for $f$ smooth of constant relative dimension $d$. $F_{*}$ is degree preserving
whereas $f^{*}$ is homogeneous of degree $d$, i.e. maps $A_{i}(\cdot)$ to 
$A_{i+d}(\cdot)$. 
\item These transformations satisfy (for such maps) the projection and
base-change formula as in the context of constructible functions $F(X)$
or of $K_{0}(var/X)$.
\item One has an associative and commutative exterior product 
\[A_{i}(X)\times A_{j}(X')\to A_{i+j}(X\times X')\]
with a unit $1_{pt}\in A_{0}(pt)$. Moreover, projective pushdown
and smooth pullback commute with exterior products:
\[(f\times f')_{*}= f_{*}\times f'_{*} \quad \text{and} \quad
 (f\times f')^{*}= f^{*}\times f'^{*} \:.\]
By these exterior products, $A_{*}(pt)$ becomes a commutative ring and
$A_{*}(X)$ a $A_{*}(pt)$-module such that $f_{*}$ and $f^{*}$ are
$A_{*}(pt)$-linear.
\item Finally one has for a line bundle $L$ over $X$
a first Chern class operation 
\[\tilde{c}_{1}(L): A_{i}(X)\to A_{i-1}(X)\]
with suitable properties, which we do not need here in full
generality. In particular, the first Chern class is nilpotent,
depends only on the isomorphism class of the line bundle, is contravariant
functorial and satisfies the projective bundle theorem
\begin{equation} \label{eq:proj}
A_{*}(\bb{P}(E)) = \sum_{i=0}^{rk\;E-1} \: \tilde{c}_{1}({\cal O}(1))^{i} \:
A_{*+i-rk\;E}(X) \tag{PB}
\end{equation} 
for $E\to X$ a vector bundle. Here ${\cal O}(1)$ is the canonical
quotient line bundle on $\bb{P}(E):=Proj(Sym^{*}{\cal E})$, with
${\cal E}$ the sheaf of sections of $E$.
\end{itemize}

A natural transformation of oriented Borel-Moore weak homology theories 
commutes by definition with these structures. By a theorem of
Levine and Morel, algebraic cobordism $\Omega_{*}$ is the universal
 example of such a theory \cite[thm.10.11,p.78]{LM}, i.e. one has a unique natural
transformation to any given oriented Borel-Moore weak homology theory.
Other examples are given by our homology theories $H_{*}(X)$
(Chow groups or Borel-Moore homology) with the usual
first Chern class operations and by $G_{0}(X)\otimes
\bb{Z}[y,y^{-1}]$. Here $y$ is of degree one, and projective pushdown 
or exterior product is defined by linear extension of the corresponding
operation for $G_{0}(X)$. But for the pullback  under a smooth morphism
$f: X'\to X$ of constant relative dimension $d$ one multiplies in addition with
$y^{d}$ to get a transformation of degree $d$:
\[f^{*}: G_{0}(X)\otimes \bb{Z}[y,y^{-1}]\to G_{0}(X')\otimes \bb{Z}[y,y^{-1}]
\:,\: [{\cal F}]y^{i}\mapsto [f^{*}{\cal F}]y^{i+d} \:.\]
The first Chern class transformation $\tilde{c}_{1}(L)$ in this example is
given by multiplication with $(1-[{\cal L}^{*}])y^{-1}$.\\

To each oriented Borel-Moore weak homology theory $A_{*}$ one can associate
a formal group law in one variable $(A_{*}(pt),F_{A}(u,v))$ \cite[sec.10.4,p.74]{LM}
given by a power series 
\[F_{A}(u,v)= \sum _{i,j=0}^{\infty}\: a_{i,j}u^{i}v^{j} \in
A_{*}(pt)[[u,v]] \quad \text{with $a_{i.j}\in A_{i+j-1}(pt)$,}\]
which describes the first Chern class operations for tensor
products of line bundles. In particular one gets an induced
graded ring homomorphism
\[\Phi_{A}: L_{*}\to A_{*}(pt)\]
from the universal Lazard ring $L_{*}$ classifying such formal group laws.
One of the main results of Levine-Morel says that $\Phi_{\Omega}$
is an isomorphism \cite[thm.12.8,p.97]{LM}. The formal group law for homology $H_{*}(X)$ is given by
the additive formal group law
\begin{equation} \label{eq:fora}
F_{H}:=F_{a}: F_{a}(u,v)=u+v\:, 
\end{equation}
and for $G_{0}(X)\otimes \bb{Z}[y,y^{-1}]$ 
it is given be the multiplicative (periodic) formal group law
\begin{equation} \label{eq:form}
F_{G}:=F_{m}: F_{m}(u,v)=u+v-y^{-1}uv \:. 
\end{equation}
Here 'periodic' stands for the fact that $y^{-1}$ is invertible.
Consider the graded ring homomorphism
\[L_{*} \to \bb{Z}[y,y^{-1}]\]
which classifies the multiplicative (periodic) formal group law (\ref{eq:form})
on $\bb{Z}[y,y^{-1}]$. Then one can form the 
oriented Borel-Moore weak homology theory
\[\Omega_{*}(X)\otimes_{\bb{L}*} \bb{Z}[y,y^{-1}] \:,\]
which is the universal such theory with 
the multiplicative (periodic) formal group law (\ref{eq:form}) \cite[rem.10.12,p.79]{LM}.
One gets in particular a natural transformation
\[\omega_{*}: \Omega_{*}(X)\otimes_{\bb{L}*} \bb{Z}[y,y^{-1}] 
\to G_{0}(X)\otimes \bb{Z}[y,y^{-1}]\:,\]
which turns out to be an isomorphism for
$X$ smooth and pure dimensional \cite[cor.11.11,p.93]{LM}. \\

This information is enough to lift our motivic Chern class
transformation $mC_{*}$ for a fixed $X$ to a unique group homomorphism
\[mC'_{*}: K_{0}(var/X) \to 
\Omega_{*}(X)\otimes_{\bb{L}*} \bb{Z}[y,y^{-1}] \:,\]
satisfying for a projective morphism  $f: X'\to X$ with $X'$ smooth
and pure-dimensional the normalization condition
\begin{equation} \label{eq:normC'}
mC'_{*}([f:X'\to X])=f_{*}(\omega_{*}^{-1}(mC_{*}([id_{X'}]))\in 
\Omega_{*}(X)\otimes_{\bb{L}*} \bb{Z}[y,y^{-1}] \:.
\end{equation}

Indeed, this defines $mC'_{*}$ on $M^{+}_{*}(X)$,
and it satisfies the blow-up relation (bl) in the
description of $K_{0}(var/X)$ given by Looijenga-Bittner \cite{Bi},
because this is first the case in $\Omega_{*}(X')\otimes_{\bb{L}*}
\bb{Z}[y,y^{-1}]$ by the isomorphism $\omega_{*}$ for smooth pure-dimensional
spaces. Note that $\omega_{*}$ commutes with projective pushdown
so that this case follows from the existence of the motivic
Chern class transformation $mC_{*}$. 
Applying the group homomorphism $f_{*}$ gives us then the relation
(bl) also in $\Omega_{*}(X)\otimes_{\bb{L}*} \bb{Z}[y,y^{-1}]$.
By functoriality of pushdown it is also clear that $mC'_{*}$ commutes with
projective pushdown! Similarly, corollary \ref{cor:mC} implies the
corresponding properties 1. and 2. also for $mC'_{*}$.\\

We now explain a more intrinsic description of the normalization condition
for the motivic Cherm class transformation $mC'_{*}$,
which allows us to formulate also the corresponding
Verdier Riemann-Roch formula for $mC'_{*}$. For this we use the
description of higher Chern classes in 
\[G_{0}(X)\otimes \bb{Z}[y,y^{-1}] \simeq 
K_{0}(X)\otimes \bb{Z}[y,y^{-1}] \]
for $X$ smooth and pure-dimensional in terms of the $\gamma$-operation
\cite{FL,SGA6}.\\

Let us first recall that higher Chern classes for vector bundles
are defined in the context of an oriented Borel-Moore weak homology theory
$A_{*}$ following the method of Grothendieck by using the projective bundle
theorem (PB). For a vector bundle $E$ over $X$ there exist unique homomorphisms
\[\tilde{c}_{i}(E): A_{*}(X)\to A_{*-i}(X) \quad i=0,\dots,rk\;E\:,\]
with $\tilde{c}_{0}(E)=1$, and satisfying the equation (in the endomorphism
ring of $A_{*}(\bb{P}(E))$):
\begin{equation} \label{eq;higherc}
\sum_{i=0}^{rk\;E}\: (-1)^{i}\cdot \tilde{c}_{i}(E)\cdot
\tilde{c}_{1}({\cal O}(1))^{rk\;E-i} \: = 0 \:.
\end{equation}
Here we omit as usual in (\ref{eq;higherc}) the pullback of 
$\tilde{c}_{i}(E)$ to $A_{*}(\bb{P}(E))$. By definition
\[\tilde{c}_{i}(E):= 0 \quad \text{for $i> rk\;E$.}\]

Note that a transformation of oriented Borel-Moore weak homology theories
(like $\omega_{*}$) is compatible with these higher Chern class operations.
Consider the total $\gamma$-class transformation \cite{FL,SGA6}
\[\gamma_{y}:=\lambda_{y/1-y}=:\sum_{i=0}^{\infty}\: \gamma^{i}y^{i}:
K^{0}(X)\to K^{0}(X)\otimes\bb{Z}[[y]] \:.\] Then one has for a vector bundle
$E$ over $X$ the well known relation  (\cite[eq.(1.13.1),p.379]{SGA6})
\begin{equation} \label{eq:gamma}
\sum_{i=0}^{rk\;E}\: \gamma^{i}(E^{*}- rk\;E)\cdot
(1-{\cal O}(1)^{*})^{rk\;E-i} \:= 0 
\end{equation}
in $K_{0}(\bb{P}(E))$. And this implies the following formula
for the higher Chern classes of $E$ in the 
oriented Borel-Moore weak homology theory $G_{0}(X)\otimes \bb{Z}[y,y^{-1}]$:
\begin{equation} \label{eq:higherc2}
\tilde{c}_{i}(E) = (-1)^{i}\cdot \gamma^{i}(E^{*}- rk\;E)\cdot y^{-i}
\end{equation} 
as a multiplication operator.\\

\begin{rem} \label{rem:gamma}
Note that this definition of higher Chern classes fits with the usual
definition of higher Chern classes on $K^{0}(X)$ with values in the
graded ring $gr_{\gamma}^{*}K^{0}(X)$ associated to the $\gamma$-filtration on
$K^{0}(X)$ (\cite{FL,SGA6}):
\[gr^{i}_{\gamma}\gamma^{i}(E- rk\;E)=:c_{i}(E) = (-1)^{i}c_{i}(E^{*})
= gr^{i}_{\gamma}(\tilde{c}_{i}(E)(y=1)) \:.\]
\end{rem}

Since the normalization condition for our motivic Chern class transformation
$mC$ was given in terms of the total $\lambda$-class $\lambda_{y}T^{*}X$ of the
cotangent bundle, we also have to express this $\lambda$-class 
in terms of the corresponding $\gamma$-classes and higher Chern classes.
Here one can use the relation (\cite[prop.1.1(a),p.48]{FL})
\[\sum_{i=0}^{rk\;E}\: \gamma^{rk\;E -i}(E^{*}-rk\;E) \cdot (1+y)^{i} =
\sum_{i=0}^{rk\;E}\: \lambda^{rk\;E-i}(E^{*}) \cdot y^{i} \:,\]
and  gets 
\[\lambda^{i}(E^{*}) = \sum_{j=0}^{i}\: \gamma^{j}(E^{*}-rk\;E)\cdot 
{rk\;E-j \choose rk\;E-i} \:.\]
In terms of higher Chern classes this gives
\[\lambda^{i}(E^{*}) = \sum_{j=0}^{i}\: \tilde{c}_{j}(E)\cdot 
{rk\;E-j \choose rk\;E-i}\cdot (-y)^{j} \:.\]
Altogether we get the following reformulation of the normalization condition
for our motivic Chern class transformation $mC_{*}$:
\begin{equation} \label{eq:reform}
\begin{split}
mC_{*}([id_{X}]) &= \left( \sum_{i=0}^{d} \: \left(\sum_{j=0}^{i}\:
\tilde{c}_{j}(TX)\cdot  {d-j \choose d-i}\cdot (-y)^{j} \right) \cdot
y^{i-d} \right)\cdot [X] \\
&= \left( \sum_{j=0}^{d} \: \left(\sum_{i=j}^{d}\: (-1)^{j}\cdot
{d-j \choose d-i}\cdot y^{i+j-d} \right) \cdot \tilde{c}_{j}(TX)
\right)\cdot [X] \\
&=: \tilde{\lambda}_{y}(T^{*}X) \cdot [X]
\in G_{0}(X)\otimes \bb{Z}[y,y^{-1}] 
\end{split} \end{equation}
for $X$ smooth and pure $d$-dimensional.
Here $[X]=k^{*}1_{pt}= [{\cal O}_{X}]y^{d} \in G_{0}(X)y^{d}$
is the 'fundamental class' of $X$ (with $k: X\to pt$ a constant map)
with respect to the  oriented Borel-Moore weak homology theory
$G_{0}(X)\otimes \bb{Z}[y,y^{-1}]$.\\

\begin{rem} \label{rem:gamma2}
We see in particular, that the normalization condition for the motivic
Chern class transformation $mC_{*}$
is not (!) of the simple form
\[mC_{*}([id_{X}]) = \left( \sum_{j=0}^{d} \: \tilde{c}_{j}(TX) \right)
 \cdot [X]\]
as it is for the Chern-Schwartz-MacPherson class transformation $c_{*}$.
We have to weight $\tilde{c}_{j}(TX)$ $(j=0,\dots,d)$ by the factor
\[\sum_{i=j}^{d}\: (-1)^{j}\cdot {d-j \choose d-i}\cdot y^{i+j-d}
\in \bb{Z}[y,y^{-1}] \:.\] 
This is related to the fact that the formal group law $F_{G}$
is the multiplicative, and not the additive formal group law.
So we have to distinguish between $\lambda$-classes and $\gamma$-classes.
Of course, this is not necessary for the usual
definition of higher Chern classes on $K^{0}(X)$ with values in the
graded ring $gr_{\gamma}^{*}K^{0}(X)$ associated to the $\gamma$-filtration on
$K^{0}(X)$:
\[c_{i}(E) = (-1)^{i}c_{i}(E^{*})
= (-1)^{i}gr^{i}_{\gamma}\gamma^{i}(E^{*}- rk\;E)
= (-1)^{i}gr^{i}_{\gamma}\lambda^{i}(E^{*}) \:.\]
\end{rem}

The right hand side of equation (\ref{eq:reform}) makes perfectly sense in 
\[\Omega_{*}(X)\otimes_{\bb{L}*} \bb{Z}[y,y^{-1}]\] 
and can be
used as a definition of the total ${\tilde{\lambda}_{y}}$-class of the
cotangent bundle (or any other vector bundle with $d:= rk\;E$) in such a way
that \[\omega_{*}(\tilde{\lambda}_{y}(T^{*}X) \cdot [X])=
\tilde{\lambda}_{y}(T^{*}X) \cdot [X] = 
\lambda_{y}([T^{*}X])\cap [{\cal O}_{X}]\:.\]
So we finally get the following

\begin{thm} \label{thm:mC'}
There exists a unique group homomorphism $mC'_{*}$ commuting with pushdown
for projective maps:
\[mC'_{*}: K_{0}(var/X)\to \Omega_{*}(X)\otimes_{\bb{L}*} \bb{Z}[y,y^{-1}]
\:,\] satisfying the normalization condition
\[mC'_{*}([id_{X}])= 
\tilde{\lambda}_{y}(T^{*}X)\cdot [X]\]
for $X$ smooth and pure-dimensional. $\quad \Box$
\end{thm}

\begin{cor} \label{cor:mC'}
\begin{enumerate}
\item $mC'_{*}$ is filtration preserving, if 
$\Omega_{*}(X)\otimes_{\bb{L}*} \bb{Z}[y,y^{-1}]$
has the filtration coming from its grading. 
\item $mC'_{*}$ commutes with exterior products.
\item One has the following Verdier Riemann-Roch formula for $f:X'\to X$
a smooth morphism of constant relative dimension:
\[\tilde{\lambda}_{y}(T^{*}_{f}) \cdot f^{*}mC'_{*}([Z\to X]) 
= mC'_{*}f^{*}([Z\to X]) \:.\]
In particular $mC'_{*}$ commutes with pullback under \'{e}tale morphisms.
$\quad \Box$
\end{enumerate} \end{cor}

$ $\\
Jean-Paul Brasselet\\
Institute de Math\'{e}matiques de Luminy\\
UPR 9016-CNRS\\
Campus de Luminy - Case 907\\
13288 Marseille Cedex 9, France\\
E-mail: jpb@iml.univ-mrs.fr\\

$ $\\
J\"{o}rg Sch\"{u}rmann\\
Westf. Wilhelms-Universit\"{a}t\\
SFB 478
"Geometrische Strukturen in der Mathematik" \\
Hittorfstr.27\\
48149 M\"{u}nster, Germany\\
E-mail: jschuerm@math.uni-muenster.de\\

$ $\\
Shoji Yokura\\
Department of Mathematics and Computer Science\\
Faculty of Science, University of Kagoshima\\
21-35 Korimoto 1-Chome\\
Kagoshima 890-0065, Japan\\
E-mail: yokura@sci.kagoshima-u.ac.jp

\end{document}